\documentclass{amsart}
\usepackage{amsmath}%
\usepackage{amsfonts}%
\usepackage{amssymb}%
\usepackage{graphicx}

\usepackage[latin1]{inputenc}
\usepackage[T1]{fontenc}

\usepackage{graphicx}
\usepackage[cmtip,arrow]{xy}
\usepackage{pb-diagram,pb-xy}

\newcommand{\field}[1]{\mathbb{#1}}
\newcommand{\german}[1]{\mathfrak{#1}}
\newcommand{\call}[1]{\mathcal{#1}}

\newcommand{\C}{\field{C}}

\newcommand{\Z}{\field{Z}}

\newcommand{\g}{\german{g}}
\newcommand{\h}{\german{h}}

\newcommand{\p}{\german{p}}
\newcommand{\q}{\german{q}}

\newcommand{\bi}{\german{b}}
\newcommand{\el}{\german{l}}
\newcommand{\uu}{\german{u}}

\newcommand{\A}{\call{A}}
\newcommand{\B}{\call{B}}
\newcommand{\Ci}{\call{C}}
\newcommand{\D}{\call{D}}
\newcommand{\E}{\call{E}}
\newcommand{\F}{\call{F}}

\newcommand{\K}{\call{K}}

\newcommand{\eS}{\call{S}}
\newcommand{\T}{\call{T}}
\newcommand{\U}{\call{U}}
\newcommand{\V}{\call{V}}
\newcommand{\W}{\call{W}}
\newcommand{\Ze}{\call{Z}}

\newcommand{\qbino}[2]{\left[
\begin{array}{c}
    #1\\
    #2
\end{array}
\right]}

 \DeclareMathOperator{\Gr}{Gr}
 
\DeclareMathOperator{\Spec}{Spec}

\DeclareMathOperator{\Rap}{Rap} 
\DeclareMathOperator{\rk}{rank}

\newcommand{\e}{\epsilon}
\newcommand{\la}{\lambda}
\newcommand{\w}{\mathbf{w}}

\newtheorem{theorem}{Theorem}
\theoremstyle{plain}

\newtheorem{corollary}{Corollary}

\newtheorem{definition}{Definition}

\newtheorem{lemma}{Lemma}

\newtheorem{proposition}{Proposition}

\numberwithin{equation}{section}

\begin{document}
\title{Degree of parabolic quantum groups}
\author{Riccardo Pulcini}
\address{Dipartimento di Matematica Università ``Roma Tre'' \newline%
\indent Largo san Leonardo Murialdo 1 J%
\indent 00146 Roma Italia}%
\email{pulcini@mat.uniroma3.it}%
\urladdr{http://xoomer.virgilio.it/pulcini\_riccardo/}

\date{14/06/2006}

\begin{abstract}
We study some elementary properties of the quantum enveloping algebra associated to a parabolic subalgebra $\p$ of a semisimple Lie algebra $\g$. In particular we prove an explicit formula for the degree of this algebra, that extends the well known formula for the quantum enveloping algebra associated to $\g$ and $\bi$, where $\bi$ is a Borel subalgebra of $\g$.
\end{abstract}

\maketitle

\section{Introduction and overview}

The aim of this work is to calculate the degree of some quantum
universal enveloping algebras. Let $\g$ be a semisimple Lie algebra,
fix a Cartan subalgebra $\h\subset\g$ and a Borel subalgebra
$\h\subset\bi\subset\g$, and denote with $\Pi$ the correspondent
set of simple roots. Given $\Pi'\subset\Pi$, we associate to
$\Pi'$ the parabolic subalgebra $\p\supset\bi$. In this situation we can "quantized" our algebras, we obtain Hopf algebras over $\C\left[q,q^{-1}\right]$, namely $\U_q(\bi)\subset\U_q(\p)\subset\U_q(\g)$.

When we specialize the parameter $q$ to a primitive $l^{th}$ root
$\e$ of 1 (with some restrictions on $l$), the resulting algebras
are finite modules over their centers, and they are a finitely generated
$\C$ algebra. In particular, every irreducible representations has
finite dimension. Let us denote by $\V$ the set of irreducible
representations, Schur lemma gives us a surjective application
$$
\pi:\V\rightarrow\Spec(Z).
$$
To determine the pull back of a point in $\Spec(Z)$ is a very
difficults problem. But generically the problem becomes easier. Since
our algebras are domains, there exists a non empty Zariski open set
$V\subset\Spec(Z)$, such that $\pi|_{\pi^{-1}(V)}$ is bijective and
moreover every irreducible representation in $\pi^{-1}(V)$ has the
same dimension $d$, the degree of our algebra. The problem is to
identify $d$.

Note that, a natural candidate for $d$ exists. We will see that in
the case of $\U_\e(\p)$, we can find a natural subalgebra
$Z_0\subset Z$, which is a Hopf subalgebra of $\U_\e(\p)$. Therefore it is the coordinate ring of an algebraic group $H$. The deformation structure of $\U_\e(p)$ implies that $H$ has a Poisson structure. Let $\delta$ be the maximal dimension of the symplectic leaves, then a natural conjecture is
$d=l^{\frac{\delta}{2}}$. This is well know in several cases, for
example, $\p=\g$ and $\p=\bi$ (cf \cite{dk1} and \cite{kw}).

Our purpose has been to prove one explicit formula for $\delta$. Before describing the strategy of the proof, we explain the formula
for $\delta$. Set $\el$ the Levi factor of $\p$. Let $\W$ be the
Weyl group of the root system of $\g$, and $\W^\el\subset\W$ that
one of subsystem generated by $\Pi'$. Denote by $\w_0$ the
longest element of $\W$ and $\w^\el_0$ the longest element of
$\W^\el$. Recall that $\W$ acts on $\h$ and set $s$ as the rank of
the linear transformation $w_0-w^\el_0$ of $\h$. Then
$$
\delta=l(w_0)+l(w^\el_0)+s,
$$
where $l$ is the length function with respect to the simple
reflection. We describe now the strategy of the proof. Our instrument is the theory of quasi polynomial algebras. By a result of De Concini, Kaç and Procesi (\cite{dkp} and \cite{dkp2}), we know that in order to compute the degree of such algebras we can be reduced to the computation of the rank of a skew symmetric matrix. To use this result we construct a degeneration of our quantum algebras to a quasipolynomial algebra, $\U^t_\e$, and a family, $\U^{t,\chi}_\e$, of finitely generated algebras parameterized by $(t,\chi)\in\C\times\Spec(Z_0)$. Then we prove that $U^0_\e$ is a quasi polynomial algebra so that the theorem of De Concini, Kaç and Procesi (\cite{dkp2}) can be applied but we notice that this give us only a lower bound for the value of $d$. To get the equality we use the family $\U^{t,\chi}_\e$. The rigidity of the semisimple algebras gives us $\U^{1,\chi}_\e\cong\U^{0,\chi}_\e$ and the theory of the algebras with trace (\cite{dp}) it tells us that the degree of our quantum algebra and $\U^{1,\chi}_\e$ are equal, and that the same is true for the quasipolynomial algebra $\U^0_\e$ and $\U^{0,\chi}_\e$, then this give the desired deduction. We close this introduction with the description of the section that compose this article. In the first section we introduce the main object of this work the quantum universal enveloping algebra associated to a parabolic Lie algebra and we give some elementary properties for this algebras. The next two sections are dedicated to the proof of the formula for the degree, in particular in section \ref{dqpa} we describe the main tool of this work, the degeneration of our quantum algebra to a quasi polynomial algebra. The last section is devoted to the study of the center of the deformation, note that the actual determination of the center of the algebra $\U_\e(\p)$ remains in general an open and potentially tricky problem. However we will propose a method, inspired by work of Premet and Skryabin (\cite{premet}), to ``lift'' elements of the center of the degenerate algebra at $t=0$ to elements of the center at least over an open set of $\Spec Z_0$, and we prove that the center is deformation invariant.

\section{Quantum enveloping algebras}

We begin by recalling some classical facts about quantum enveloping algebras associated to a simple Lie algebra $\g$, and we introduce the main object of this work, the quantum enveloping algebra associated to a parabolic subalgebra $\p\subset\g$. 

Choose a Cartan subalgebra $\h\subset\g$ and a Borel subalgebra $\bi\subset\g$. Let $R\supset R^+\supset\Pi=\left\{\alpha_1,\ldots,\alpha_n\right\}$ be the root system, set of positive root and set of simple root Let $C=(c_{ij})$ be the Cartan matrix associated to this choice and $\underline{d}=(d_i)$ be the vector of relative positive integer associated to $C$ such that $(d_ic_{ij})$ is a symmetric matrix. As usual we call $\W$ the Weyl group associated to $\h$ and $\bi$ and we set $\w_0$ as the longest element of $\W$. Denote by $P$ and $Q$ the its weight and root lattices and let $w_1,\ldots,w_n$ be the fundamental weight. 

Fix $\Pi'\subset\Pi$, we call $\p$ the parabolic subalgebra associated to it. Note that if $\Pi'=\emptyset$ then $\p=\bi$ and if $\Pi'=\Pi$ then $\p=\g$. Let $\p=\el\oplus\uu$ be the Levi decomposition of $\p$, with $\el$ the Levi factor and $\uu$ the unipotent part. We call $\W^\el\subset\W$ the Weyl group of $\el$ and $\w^\el_0$ its longest element, and we have $\Pi'=\Pi^\el$.

\subsection{The simple case}

Following De Concini and Kac (\cite{dk1}), we define:
\begin{definition}\label{quantumgr}
A \emph{simply connected quantum group}\index{Quantum group}
$\U_q(\g)$ associated to the Cartan matrix $C=(c_{i,j})_{i,j=1,\ldots,n}$ is an algebra over $\C\left(q\right)$ on generators $E_{i}$, $F_{i}$ $(i=1,\cdots ,n)$, $K_{\alpha}$ $\alpha\in P$, subject to the following relations
\begin{eqnarray}
    &&\left\{\begin{array}{l}
        K_{\alpha}K_{\beta}=K_{\alpha+\beta}\\
        K_{0}=1
    \end{array}\right.\\
    \label{comm1}&&\left\{\begin{array}{l}
        K_\alpha E_{i}K_{-\alpha}=q^{\left(\alpha|\alpha_{i}\right)}E_{i}\\
        K_{\alpha}F_{i}K_{-\alpha}=q^{-\left(\alpha|\alpha_{i}\right)}F_{i}
    \end{array}\right.\\
    &&\label{comm2}\left[E_{i},F_{j}\right]=\delta_{ij}\frac{K_{\alpha_{i}}-K_{-\alpha_{i}}}{q^{d_{i}}-q^{-d_{i}}}\\
    &&\left\{\begin{array}{l}
        \sum\limits^{1-a_{ij}}_{s=0}(-1)^s\qbino{1-a_{ij}}{s}_{d_i}E^{1-a_{ij}-s}_{i}E_jE^{s}_{i}=0\mbox{ if }i\neq j\\
        \sum\limits^{1-a_{ij}}_{s=0}(-1)^s\qbino{1-a_{ij}}{s}_{d_i}F^{1-a_{ij}-s}_{i}F_jF^{s}_{i}=0\mbox{ if }i\neq j
    \end{array}\right.
\end{eqnarray}
Here is $\qbino{n}{m}_{d}=\frac{[n]_d!}{[k]_d![n-k]_d!}$ and $[n]_d=\frac{q^n-q^{-n}}{q^d-q^{-d}}$.
\end{definition}

It is well known, by the work of Lusztig (\cite{LusBook}), that
\begin{theorem}\label{Hopfstructure}
$\U_q(\g)$ has a Hopf algebra structure with comultiplication $\Delta$,
antipode $S$ and counit $\eta$ defined by:
\begin{itemize}
    \item $\left\{\begin{array}[l]{l}
        \Delta (E_{i})=E_{i}\otimes 1+K_{\alpha_{i}}\otimes E_{i}\\
        \Delta (F_{i})=F_{i}\otimes K_{-\alpha_{i}}+1\otimes F_{i}\\
        \Delta (K_{\alpha})=K_{\alpha}\otimes K_{\alpha}
    \end{array}\right.$
    \item $\left\{\begin{array}[l]{l}
        S(E_{i})=-K_{\alpha_{i}}E_{i}\\
        S(F_{i})=-F_{i}K_{\alpha_{i}}\\
        S(K_{\alpha})=K_{-\alpha}
    \end{array}\right.$
    \item $\left\{\begin{array}[l]{l}
        \eta (E_{i})=0\\
        \eta (F_{i})=0\\
        \eta (K_{\alpha})=1
    \end{array}\right.$
\end{itemize}
\end{theorem}

We denote by $\U^+$, $\U^-$ and $\U^0$ the $\C(q)$-subalgebra generated by the $E_i$, the $F_i$ and $K_\beta$ respectively. The algebras $\U^+$ and $\U^-$ are not Hopf subalgebras. On the other hand, the algebras $\U^{\geq 0}:=\U^+\U^0$ and $\U^{\leq 0}:=\U^0\U^-$ are Hopf subalgebras and we shall think to them as quantum deformation of the enveloping algebras $\U(\bi)$ and $\U(\bi^-)$. We denote them $\U_q(\bi)$ and $\U_q(\bi^-)$.

Following Lusztig (\cite{LusBook}), we define an action of the braid group $\B_\W$ (associated to $\W$). Denote by $T_i$ and $s_i$ the canonical generators of $\B_\W$ and $\W$, we define the action as an automorphism of $\U_\q(\g)$, by the formulas:
\begin{eqnarray}\label{quantumact}
    &&T_iK_\la=K_{s_i(\la)}\\
    &&T_iE_i=-F_iK_i\\
    &&T_iF_i=-K^{-1}_{i}E_i\\
    &&T_iE_j=\sum^{-c_{ij}}_{s=0}(-1)^{s-c_{ij}}q^{-sd_i}\frac{E^{-c_{ij}-s}_i}{[-c_{ij}-s]_{d_i}!}E_j\frac{E^s_i}{[s]_{d_i}!}\\
    &&T_iF_j=\sum^{-c_{ij}}_{s=0}(-1)^{s-c_{ij}}q^{sd_i}\frac{F^{s}_i}{[s]_{d_i}!}F_j\frac{F^{-c_{ij}-s}_i}{[-c_{ij}-s]_{d_i}!}.
\end{eqnarray}

We use the braid group to construct analogues of the root
vectors associated to non simple roots.

Take a reduced expression $\w_0=s_{i_1}\ldots s_{i_N}$ for
the longest element in the Weyl group $\W$. Setting
$\beta_j=s_{i_1}\cdots s_{i_{j-1}}(\alpha_j)$, we get a total order
on the set of positive root. We define the elements
$E_{\beta_j}=T_{i_1}\ldots T_{i_{j-1}}(E_{i_j})$ and
$F_{\beta_j}=T_{i_1}\ldots T_{i_{j-1}}(F_{i_j})$. Note that this
elements depend on the choice of the reduced expression.

\begin{lemma}
\begin{enumerate}
    \item $E_{\beta_j}\in \U^+$, $\forall$ $i=1\ldots N$ and the monomials $E^{k_1}_{\beta_1}\cdots E^{k_N}_{\beta_N}$ form a $\C(q)$ basis of $\U^+$
    \item $F_{\beta_j}\in \U^-$, $\forall$ $i=1\ldots N$ and the monomials $F^{k_1}_{\beta_1}\cdots F^{k_N}_{\beta_N}$ form a $\C(q)$ basis of $\U^-$
\end{enumerate}
\end{lemma}

\begin{theorem}[Poincaré-Birkoff-Witt theorem]\label{PBWthm}
The monomials
$$
E^{k_1}_{\beta_1}\cdots E^{k_N}_{\beta_N}K_\alpha F^{k_N}_{\beta_N}\cdots F^{k_1}_{\beta_1}
$$
are a $\C(q)$ basis of $\U$. In fact as vector spaces, we have the
tensor product decomposition,
$$
\U=\U^+\otimes\U^0\otimes\U^-
$$
\end{theorem}
\begin{proof}
See \cite{LusBook}.
\end{proof}

\begin{theorem}[Levendorskii-Soibelman relations]\label{LSrel}
For $i<j$ one has
\begin{enumerate}
    \item
    \begin{equation}\label{LS1}
    E_{\beta_j}E_{\beta_i}-q^{(\beta_i|\beta_j)}E_{\beta_i}E_{\beta_j}=\sum_{k\in\Z^{N}_{+}}c_kE^k
    \end{equation}
    where $c_k\in \C[q,q^{-1}]$ and $c_k\neq 0$ only when $k=\left(k_1,\ldots,k_N\right)$ is such that $k_s=0$ for $s\leq i$ and $s\geq j$, and $E^k=E^{k_1}_{\beta_1}\cdots E^{k_N}_{\beta_N}$.
    \item
    \begin{equation}
    F_{\beta_j}F_{\beta_i}-q^{-(\beta_i|\beta_j)}F_{\beta_i}F_{\beta_j}=\sum_{k\in\Z^{N}_{+}}c_kF^k
    \end{equation}
    where $c_k\in \C[q,q^{-1}]$ and $c_k\neq 0$ only when $k=\left(k_1,\cdots,k_N\right)$ is such that $k_s=0$ for $s\leq i$ and $s\geq j$, and $F^k=F^{k_N}_{\beta_N}\cdots F^{k_1}_{\beta_1}$.
\end{enumerate}
\end{theorem}
\begin{proof}
See \cite{LS}.
\end{proof}

To obtain from $\U_q(\g)$ a well defined Hopf algebra by specializing $q$
to an arbitrary non zero complex number $\e$, one can construct an
\emph{integral form} of $\U$.
\begin{definition}
An \emph{integral form} $\U_\A$ is a $\A$ subalgebra, where $\A=\C[q,q^{-1}]$, such that the natural map
$$
\U_\A\otimes_\A\C(q)\mapsto\U
$$
is an isomorphism of $\C(q)$ algebra.
We define
$$
\U_\e=\U_\A\otimes_\A\C
$$
using the homomorphism $\A\mapsto\C$ mapping $q$ to $\e$.
\end{definition}

There are two different candidates for $\U_\A$ the non restricted
and the restricted integral form, which lead to different
specializations (with markedly different representation theories)
for certain values of $\e$. We are interested in the non restricted
form, for more details one can see \cite{Poisson1}.

Introduce the elements
$$
\left[K_i;m\right]_{q_i}=\frac{K_iq^{m}_{i}-K^{-1}_{i}q^{-m}_{i}}{q_i-q^{-1}_{i}}\in\U^{0}
$$
with $m\geq 0$, where $q_{i}=q^{d_i}$.
\begin{definition}\label{quantumunita}
The algebra $\U_\A$ is the $\A$ subalgebra of $\U$ generated by the
elements $E_i$, $F_i$, $K^{\pm1}_{i}$ and
$L_i=\left[K_i;0\right]_{q_i}$, for $i=1,\ldots,n$. With the map
$\Delta$, $S$ and $\eta$ defined on the first set of generators as
in \ref{Hopfstructure} and with
\begin{eqnarray*}
    &&\Delta(L_i)=L_i\otimes K_i+K^{-1}_{i}\otimes L_i\\
    &&S(L_i)=-L_i\\
    &&\eta(L_i)=0
\end{eqnarray*}
\end{definition}

The defining relation of $\U_\A$ are as in \ref{quantumgr} replacing
\ref{comm2} by
$$
E_iF_j-F_jE_i=\delta_{ij}L_i
$$
and adding the relation
$$
(q_i-q^{-1}_{i})L_i=K_i-K^{-1}_{i}
$$

\begin{proposition}
$\U_\A$ with the previous definition is a Hopf algebra. Moreover,
$\U_\A$ is an integral form of $\U$.
\end{proposition}
\begin{proof}
See \cite{Poisson1} or \cite{dp} $§12$.
\end{proof}
\begin{proposition}\label{proUe}
If $\e^{2d_i}\neq1$ for all $i$, then
\begin{enumerate}
    \item $\U_\e$ is generated over $\C$ by the elements $E_i$, $F_i$, and $K^{\pm1}_{i}$ with defining relations obtained from those in \ref{quantumgr} by replacing $q$ by $\e$
    \item The monomials
$$
E^{k_1}_{\beta_1}\cdots E^{k_N}_{\beta_N}K_\alpha F^{k_N}_{\beta_N}\cdots F^{k_1}_{\beta_1}
$$
form a $\C$ basis of $\U_\e$.
    \item The LS relations holds in $\U_\e$.
\end{enumerate}
\end{proposition}
\begin{proof}
See \cite{dp} $§12$.
\end{proof}

\subsection{Parabolic case}

Choose a reduced expression of $w_0=s_{j_1}\ldots s_{j_k}s_{i_1}\ldots s_{i_h}$, such that $w^\el_0=s_{i_1}\ldots s_{i_h}$ is a reduced expression for $w^{\el}_{0}$. Set $\overline{w}=w_0(w^{\el}_{0})^{-1}=s_{j_1}\ldots s_{j_k}$, with
$h=|(R^\el)^+|$ and $h+k=N=|R^+|$. Define, as in the general case, 
\begin{eqnarray*}
\beta^1_t&=&\overline{w}s_{i_1}\ldots s_{i_{t-1}}\left(\alpha_{i_t}\right)\in(R^\el)^+,\\
\beta^2_t&=&s_{j_1}\ldots s_{j_{t-1}}\left(\alpha_{i_{t+k}}\right)\in R^+\setminus (R^\el)^+.
\end{eqnarray*}

Given this choice of positive root, we obtain the following $q$ analogues of the root vectors:
\begin{eqnarray*}
E_{\beta^1_t}&=&T_{\overline{w}}T_{i_1}\ldots T_{i_{t-1}}\left(E_{i_t}\right),\\
E_{\beta^2_t}&=&T_{j_1}\ldots T_{j_{t-1}}\left(E_{i_{t+k}}\right).
\end{eqnarray*}
and
\begin{eqnarray*}
F_{\beta^1_t}&=&T_{\overline{w}}T_{i_1}\ldots T_{i_{t-1}}\left(F_{i_t}\right),\\
F_{\beta^2_t}&=&T_{j_1}\ldots T_{j_{t-1}}\left(F_{i_{t+k}}\right).
\end{eqnarray*}
The PBW theorem implies that the monomials
\begin{equation}\label{pbwbasisp}
E^{s_1}_{\beta^{2}_{1}}\cdots E^{s_k}_{\beta^{2}_{k}}E^{s_{k+1}}_{\beta^{1}_{1}}\cdots E^{s_{k+h}}_{\beta^{1}_{h}}K_{\la}F^{t_{k+h}}_{\beta^{1}_{h}}\cdots F^{t_{k+1}}_{\beta^{1}_{1}}F^{t_k}_{\beta^{2}_{k}}\cdots F^{t_1}_{\beta^{2}_{1}}
\end{equation}
for $(s_1,\cdots,s_N)$, $(t_1,\ldots,t_N)\in(\Z^+)^N$ and
$\la\in\Lambda$, form a $\C(q)$-basis of $\U_q(\g)$.

The choice of the reduced expression of $w_0$ and the LS relations for $\U_q(\g)$ implies that
\begin{proposition}
For $i<j$ one has
\begin{enumerate}
    \item
    $$
    E_{\beta^1_j}E_{\beta^1_i}-q^{(\beta^1_i|\beta^1_j)}E_{\beta^1_i}E_{\beta^1_j}=\sum_{k\in\Z^{N}_{+}}c_kE_1^k
    $$
    where $c_k\in \C(q)$ and $c_k\neq 0$ only when $k=\left(s_1,\ldots,s_k\right)$ is such that $s_r=0$ for $r\leq i$ and $r\geq j$, and $E_1^k=E^{s_1}_{\beta^1_1}\ldots E^{s_k}_{\beta^1_k}$.
    \item
    $$
    E_{\beta^2_j}E_{\beta^2_i}-q^{-(\beta^2_i|\beta^2_j)}E_{\beta^2_i}E_{\beta^2_j}=\sum_{k\in\Z^{N}_{+}}c_kE_2^k
    $$
    where $c_k\in \C(q)$ and $c_k\neq 0$ only when $k=\left(t_1,\ldots,t_h\right)$ is such that $t_r=0$ for $r\leq i$ and $r\geq j$, and $E_2^k=E^{t_1}_{\beta^2_1}\ldots E^{t_h}_{\beta^2_h}$.
\end{enumerate}
The same statement holds for $F_{\beta^{1}_{i}}$ and
$F_{\beta^{2}_{i}}$.
\end{proposition}
Let $\Pi^\el$ be the simple root associated to the Levi factor $\el$. The definition of the braid group action implies:
\begin{proposition}
\begin{enumerate}
    \item If $i\in \Pi^\el$, then $E_i=E_{\beta^{1}_{s}}$ and $F_i=F_{\beta^{1}_{s}}$, for some $s\in\left\{1,\ldots,h\right\}$.
		\item If $i\in\Pi\setminus\Pi^\el$, then $E_i=E_{\beta^{2}_{s}}$ and   $F_i=F_{\beta^{2}_{s}}$ for some $s\in\left\{1,\ldots,k\right\}$.
\end{enumerate}
\end{proposition}
\begin{definition}
The simple connected quantum group associated to $\p$, or
\emph{parabolic quantum group}, is
the $\C(q)$ subalgebra of $\U_q(\g)$ generated by
$$
\U_q(\p)=\langle E_{\beta^1_i},K_\la,F_{\beta_j}\rangle
$$
for $i=1,\ldots,h$, $j=1\ldots N$ and $\la\in \Lambda$.
\end{definition}
\begin{definition}
\begin{enumerate}
    \item The quantum Levi factor of $\U_q(\p)$ is the subalgebra generated by
$$
\U_q(\el)=\langle E_{\beta^1_i},K_\la,F_{\beta^1_i}\rangle
$$
for $i=1,\ldots,h$, and $\la\in\Lambda$.
    \item The quantum unipotent part of $\U(\p)$ is the subalgebra generated by
    $$
    \U^{\overline{w}}=\langle F_{\beta^2_s}\rangle
    $$
    with $s=1\ldots h$
\end{enumerate}
\end{definition}
Set $\U_q^+(\p)=\U_q^+(\el)=\langle E_i\rangle_{i\in\Pi^\el}$, $\U_q^-(\p)=\langle F_i\rangle_{i\in\Pi}$, $\U_q^-(\el)=\langle F_i\rangle_{i\in\Pi^\el}$ and $\U_q^0(\p)=\U_q^0(\el)=\langle K_\la\rangle_{\la\in\Lambda}$. We have:
\begin{proposition}
The definition of $\U_q(\p)$ and $\U_q(\el)$ is independent of the
choice of the reduced expression of $w^\el_0$ and $w_0$.
\end{proposition}
\begin{proof}
Follows immediately from proposition 9.3 in \cite{dp}.
\end{proof}
We can  easely see taht the PBW theorem and the LS relations holds in $\U_q(\p)$ and $\U_q(\el)$, which is an immediately consequence of \ref{pbwbasisp}.
\begin{proposition}
Set $m=\rk\el=\#|\Pi^\el|$. The algebra $\U(\p)$ is generated by
$E_i$, $F_j$ $K_\la$, with $i=1,\ldots,m$, $j=1,\ldots,n$ and
$\la\in\Lambda$, subject to the following relations:
\begin{eqnarray}
    &&\left\{\begin{array}{l}
        K_{\alpha}K_{\beta}=K_{\alpha+\beta}\\
        K_{0}=1
    \end{array}\right.\\
    &&\left\{\begin{array}{l}
        K_{\alpha}E_{i}K_{-\alpha}=q^{\left(\alpha|\alpha_{i}\right)}E_{i}\\
        K_{\alpha}F_{j}K_{-\alpha}=q^{-\left(\alpha|\alpha_{j}\right)}F_{j}
    \end{array}\right.\\
    &&\left[E_{i},F_{j}\right]=\delta_{ij}\frac{K_{\alpha_{i}}-K_{-\alpha_{i}}}{q^{d_{i}}-q^{-d_{i}}}\\
    &&\left\{\begin{array}{l}
        \sum\limits^{1-a_{ij}}_{s=0}(-1)^s\qbino{1-a_{ij}}{s}_{d_i}E^{1-a_{ij}-s}_{i}E_jE^{s}_{i}=0\mbox{ if }i\neq j\\
        \sum\limits^{1-a_{ij}}_{s=0}(-1)^s\qbino{1-a_{ij}}{s}_{d_i}F^{1-a_{ij}-s}_{i}F_jF^{s}_{i}=0\mbox{ if }i\neq j.
    \end{array}\right.
\end{eqnarray}
\end{proposition}
\begin{proof}
Follows from PBW theorem and the LS relations.
\end{proof}
We state now some easy properties of $\U_q(\p)$:
\begin{lemma}
\begin{enumerate}
\item The multiplication map
$$
\U^+(\el)\otimes\U^0(\el)\otimes\U^-(\el)\rightarrow\U(\el)
$$
is an isomorphism of vector spaces.
\item The multiplication map
$$
\U(\el)\otimes \U^{\overline{w}}\stackrel{m}{\longrightarrow}\U(p)
$$
defined by $m(x,u)=xu$ for every $x\in\U(\el)$ and
$u\in\U^{\overline{w}}$, is an isomorphism of vector spaces.
\item The map $\mu:\U(\p)\rightarrow\U(\el)$ defined by
\begin{equation*}
\begin{split}
\mu\left(E^{s_{1}}_{\beta^{1}_{1}}\cdots E^{s_{h}}_{\beta^{1}_{h}}K_{\la}F^{t_{k+h}}_{\beta^{1}_{h}}\cdots F^{t_{k+1}}_{\beta^{1}_{1}}F^{t_k}_{\beta^{2}_{k}}\cdots F^{t_1}_{\beta^{2}_{1}}\right)\\
=\left\{\begin{array}{ll}
0 & \mbox{if }t_{k+i}\neq0\mbox{ for some }i=1,\ldots,h,\\
E^{s_{1}}_{\beta^{1}_{1}}\cdots E^{s_{h}}_{\beta^{1}_{h}}K_{\la}F^{t_{k+h}}_{\beta^{1}_{h}}\cdots F^{t_{k+1}}_{\beta^{2}_{1}} & \mbox{if }t_{k+i}=0\mbox{ for all }i=1,\ldots,h,
\end{array}\right.
\end{split}
\end{equation*}
is an homomorphism of algebras.
\item $\U(\p)$ and $\U(\el)$ are Hopf subalgebras of $\U$.
\end{enumerate}
\end{lemma}
\begin{proof}
Follows immediately from the definition.
\end{proof}

Let $\A=\C[q,q^{-1}]$, and $\U_\A$ the integral form of $\U_q(\g)$ defined
in definition \ref{quantumunita}. Like in the general case, we define
$\U_\A(\p)$, has the subalgebra generated by $E_{\beta^{1}_{i}}$,
$F_{\beta^{1}_{i}}$, $F_{\beta^{2}_{s}}$, $K^{\pm1}_{j}$ and $L_j$,
with $i=1,\ldots,h$, $s=1,\ldots,k$ and $j=1,\ldots,n$.

\begin{definition}
Let $\e\in\C$, we define
$$
\U_\e(\p)=\U_\A(\p)\otimes_\A\C
$$
using the homomorphism $\A\rightarrow\C$ mapping $q\rightarrow\e$
\end{definition}

Let $\e\in\C$ such that $\e^{2d_i}\neq1$ for all $i$, then
\begin{proposition}
$\U_\e(\p)\subset\U_\e(\g)$. Moreover $\U_\e(\p)$ is generated by
$E_{\beta^1_i}$, $F_{\beta_s}$ and $K^{\pm1}_{j}$, for
$i=1,\ldots,h$, $s=1,\ldots,N$ and $j=1,\ldots,n$.
\end{proposition}
\begin{proof}
The claim is a consequence of the definition of $\U_\A(\p)$.
\end{proof}
\begin{proposition}\label{LSp}
The PBW theorem and the LS relations holds for $\U_\e(\p)$
\end{proposition}
\begin{proof}
The claim is a consequence of the PBW theorem and LS relations
for $\U_\e(\g)$ and the choice of the decomposition of the reduced
expression of $w_0$.
\end{proof}

\subsection{Some observations on the center of $\U_\e(\p)$}

The aim of this section is to extend some properties of the center of $\U_\e(\g)$ at the center of $\U_\e(\p)$.

\begin{proposition}
For $i=1,\ldots,k$, $s=1,\ldots,h$ and $j=1,\ldots,n$,
$E^{l}_{\beta^{1}_{i}}$, $F^{l}_{\beta^{1}_{i}}$,
$F^{l}_{\beta^{2}_{s}}$ and $K^{\pm l}_{j}$ lie in the center of
$\U_\e(\p)$
\end{proposition}
\begin{proof}
It is well known that these elements lie in the center of $\U_\e$ (cf. \cite{dp}), but they also lie in $\U_\e(\p)$,
hence the claim.
\end{proof}

For $\alpha\in (R^\el)^+$, $\beta\in R^+$ and $\la\in Q$, define
$e_\alpha=E^{l}_{\alpha}$, $f_\beta=F^{l}_{\beta}$,
$k^{\pm1}_{\la}=K^{\pm l}_{\la}$. Let $Z_{0}(\p)$ be the subalgebra
generated by the $e_\alpha$, $f_\beta$ and $k^{\pm1}_{i}$.
\begin{proposition}\label{center2}
Let $Z^{0}_{0}$, $Z^{+}_{0}$ and $Z^{-}_{0}$ be the subalgebra
generated by $k^{\pm1}_{i}$, $e_\alpha$ and $f_\beta$ respectively.
\begin{enumerate}
    \item $Z^{\pm}_{0}\subset\U^{\pm}_{\e}(\p)$
    \item Multiplication defines an isomorphism of algebras
    $$
    Z^{-}_{0}\otimes Z^{0}_{0}\otimes Z^{+}_{0}\rightarrow Z_{0}(\p)
    $$
    \item\label{nothZ0} $Z^{0}_{0}$ is the algebra of Laurent polynomial in the $k_i$, and $Z^{+}_{0}$ and $Z^{-}_{0}$ are polynomial algebra with generators $e_\alpha$ and $f_\beta$ respectively.
    \item $\U_\e(\p)$ is a free $Z^{0}_{\e}(\p)$ module with basis the set of monomial
    $$
    E^{s_{1}}_{\beta^{1}_{1}}\cdots E^{s_{h}}_{\beta^{1}_{h}}K^{r_1}_{1}\cdots K^{r_n}_{n}F^{t_{k+h}}_{\beta^{1}_{h}}\cdots F^{t_{k+1}}_{\beta^{1}_{1}}F^{t_k}_{\beta^{2}_{k}}\cdots F^{t_1}_{\beta^{2}_{1}}
    $$
    for which $0\leq s_j,t_i,r_v<l$
\end{enumerate}
\end{proposition}
\begin{proof}
By definition of $\U^+(\p)$, we have $e_\alpha\in\U^+(\p)$, since
$\U^+(\p)$ is a subalgebra (i) follows. (ii) and (iii) are easy
corollaries of the definitions and of the PBW theorem. (iv) follows from the PBW theorem for $\U(\p)$.
\end{proof}

The previous proposition shows that $\U_\e(\p)$ is a finite
$Z_{0}(\p)$ module. Since $Z_0$ is clearly Noetherian, from
\ref{nothZ0}, it follows that $Z_{\e}(\p)\subset \U_\e(\p)$ is a
finite $Z_{0}(\p)$ module, and hence integral over $Z_{0}(\p)$. By
the Hilbert basis theorem $Z_{\e}(\p)$ is a finitely generated
algebra. Thus the affine schemes $\Spec(Z_\e(\p))$ and
$\Spec(Z_{0}(\p))$ are algebraic varieties. Note that $\Spec(Z_{0})$
is isomorphic to $\C^{N}\times\C^{l(h)}\times(\C^{*})^n$. Moreover
the inclusion $Z_{0}(\p)\hookrightarrow Z_\e(\p)$ induces a
projection $\tau:\Spec(Z_\e(\p))\rightarrow\Spec(Z_{0}(\p))$, and we
have
\begin{proposition}
$\Spec(Z_\e(\p))$ is an affine variety and $\tau$ is a finite
surjective map.
\end{proposition}
\begin{proof}
Follows from the Cohen-Seidenberg theorem (\cite{serre} ch. III).
\end{proof}

We conclude this section by discussing the relation between the
center and the Hopf algebra structure of $\U_\e(\p)$.

\begin{proposition}
\begin{enumerate}
    \item $Z_{0}(\p)$ is a Hopf subalgebra of $\U_\e(\p)$.
    \item $Z_{0}(\p)$ is a Hopf subalgebra of $Z_{0}$.
    \item $Z_{0}(\p)=Z_{0}\cap\U_\e(\p)$.
\end{enumerate}
\end{proposition}
\begin{proof}
It follows directly from the given definitions.
\end{proof}
The fact that $Z_0(\p)$ is an Hopf algebra tells us that
$\Spec(Z_{0}(\p))$ is an algebraic group . Moreover, the inclusion
$Z_0(\p)\hookrightarrow Z_0$ being an inclusion of Hopf algebras,
induces a group homomorphism,
$$
\Spec(Z_0)\rightarrow\Spec(Z_0(\p)).
$$
Let us recall that in \cite{dkp} the authors prove that the center $Z_0$ has the following form:
$$
\Spec(Z_0)=\left\{(a,b):\in B^-\times B^+:\pi^-(a)\pi^+(b)=1\right\}
$$
where if we denote by $G$ the connected simply connected Lie group associated to $\g$, then $B^\pm$ are the borel subgroups of $G$ associated to $\bi^\pm$, $H$ is the toral subgroup associated to $\h$ and $\pi^\pm:B^\pm\rightarrow H$ is the canonical map. From this and, the explicit description of the subalgebra $Z_0(\p)\subset Z_0$, we get
$$
\Spec(Z_0(\p))=\left\{(a,b):\in B^-_L\times B^+:\pi^-(a)\pi^+(b)=1\right\}
$$
where $L\subset G$ is the connected subgroup of $G$ such that $Lie(L)=\el$, and $B^-_L=B^-\cap L$.

\section{Degeneration to a quasi polynomial algebra}\label{dqpa}

\subsection{The case $\p=\g$.}\label{gcaso}

\begin{definition}
Let $t\in \C$, we define $\U^{t}_{\e}$ the algebra over $\C$ on
generators $E_i$, $F_i$, $L_i$ and $K^{\pm}_{i}$, for
$i=1,\ldots,n$, subject to the following relations:
\begin{eqnarray}\label{quantumreldef}
    &&\left\{\begin{array}{l}
        K^{\pm1}_{i}K^{\pm1}_{j}=K^{\pm1}_jK^{\pm1}_{i}\\
        K_{i}K^{-1}_{i}=1
    \end{array}\right.\\
    &&\left\{\begin{array}{l}
        K_iE_{j}K^{-1}_{i}=\e^{a_{ij}}E_{j}\\
        K_iF_{j}K^{-1}_{i}=\e^{-a_{ij}}F_{j}
    \end{array}\right.\\
    &&\left\{\begin{array}{l}
        \left[E_{i},F_{j}\right]=t\delta_{ij}L_i\\
        \left(ad_{\sigma_{-\alpha_{i}}}E_{i}\right)^{1-a_{ij}}E_{j}=0\\
        \left(ad_{\sigma_{-\alpha_{i}}}F_{i}\right)^{1-a_{ij}}F_{j}=0
    \end{array}\right.\\
    &&\left\{\begin{array}{l}
        \left(\e^{d_i}-\e^{-d_i}\right)L_i=t\left(K_i-K^{-1}_{i}\right)\\
        \left[L_i,E_j\right]=t\frac{\e^{a_{ij}}-1}{\e^{d_i}-\e^{-di}}\left(E_jK_i+K^{-1}_{i}E_j\right)\\
        \left[L_i,F_j\right]=t\frac{\e^{-a_{ij}}-1}{\e^{d_i}-\e^{-di}}\left(F_jK_i+K^{-1}_{i}F_j\right)
    \end{array}\right.
\end{eqnarray}
\end{definition}
Let $0\neq\lambda\in \C$, define
\begin{equation}\label{vartheta}
\vartheta_\la(E_i)=\frac{1}{\la}E_i,\mbox{ }\vartheta_\la(F_i)=\frac{1}{\la}F_i,\mbox{ }\vartheta_\la(L_i)=\frac{1}{\la}L_i,\mbox{ }\vartheta_\la(K^{\pm1}_i)=K^{\pm1}_i,
\end{equation}
for $i=1,\ldots,n$.
\begin{proposition}
For any $0\neq \la\in\C$, $\vartheta_\la$ is an isomorphism of
algebra between $\U^{t}_{\e}$ and $\U^{\la t}_{\e}$. In particular if $t\neq0$ then $\U^t_\e\cong\U_\e(\g)$.
\end{proposition}
\begin{proof}
Simple verification of the properties.
\end{proof}
Set $\eS_\e:=\U^{t=0}_\e$, we want to construct an explicit
realization of it. Let $\D=\U_{\e}(\bi_+)\otimes \U_{\e}(\bi_-)$ and
define the map
$$
\Sigma:\eS_\e\rightarrow \D
$$
by $\Sigma(E_i)=\E_i:=E_i\otimes 1$, $\Sigma(F_i)=\F_i=1\otimes F_i$, and $\Sigma(K^{\pm 1}_i)=\K^{\pm 1}_{i}:=K^{\pm 1}_i\otimes K^{\pm 1}_i$ for $i=1,\ldots,n$.
\begin{lemma}\label{sigma}
$\Sigma$ is a well defined map.
\end{lemma}
\begin{proof}
We must verify that the image of $E_i$, $F_i$ and $K_i$ satisfy the
relation \ref{quantumreldef} for $t=0$. This is a simple matter of bookkeeping.
\end{proof}

Note that $\Sigma$ is injective, then we can identify $\eS_\e$ with
the subalgebra of $\D$ generated by $\E_i$, $\F_i$ and $\K_i$, for
$i=1,\ldots,n$. We define now the analogues of the root vectors for
$\eS_\e$:
\begin{definition}
For all $i=1,\ldots,N$, let
\begin{enumerate}
    \item $\E_{\beta_i}:=E_{\beta_i}\otimes 1\in \eS_\e$
    \item $\F_{\beta_i}:=1\otimes F_{\beta_i}\in \eS_\e$
\end{enumerate}
\end{definition}
As a consequence of this we get a PBW theorem for $\eS_\e$.
\begin{proposition}\label{PBWbasisS}
The monomials
$$
\E^{k_1}_{\beta_1}\ldots \E^{k_N}_{\beta_N}\K^{s_1}_{1}\ldots \K^{s_n}_{n}\F^{h_1}_{\beta_N}\ldots \F^{k_1}_{\beta_1}
$$
for $(k_1,\ldots,k_N)$, $(h_1,\ldots,h_N)\in(\Z^{+})^{N}$ and
$(s_1,\ldots,s_n)\in\Z^{n}$, form a $\C$ basis of $\eS_\e$. Moreover
$$
\eS_\e=\eS^-_\e\otimes \eS^0_\e\otimes\eS^+_\e
$$
where $\eS^+_\e$ (resp. $\eS^-_\e$ and $\eS^0_\e$) is the subalgebra
generated by $\E_{\beta_i}$ (resp. $\F_{\beta_i}$ and $\K_i$).
\end{proposition}
\begin{proof}
This follows from the injectivity of $\Sigma $ and PBW theorem for $\U_\e(\g)$
\end{proof}

It is clear that $\E_{\beta_i}$ is also the image of the element
$E_{\beta_i}\in\U^t_\e$, where the $E_{\beta_i}$ are non commutative
polynomials in the $E_i$'s by Lusztig procedure (\cite{LusBook}).
The same thing is true for $\F_{\beta_i}$ and $F_{\beta_i}$.

It is also clear that the LS relations hold for $\E_{\beta_i}$ and the $\F_{\beta_i}$ (instead of the $E_{\beta_i}$ and the $F_{\beta_i}$) inside $\eS_\e$. 

\begin{theorem}\label{defor}
$\eS_\e=\U^{t=0}_{\e}$ is a twisted derivation algebra.
\end{theorem}
\begin{proof}
Define $\U^{0}=\C[\E_{\beta_1},\F_{\beta_N}]\subset \eS_\e$, then we can define
$$
\U^{i}=\U^{i-1}_{\sigma,D}\left[\E_{\beta_{i}},\F_{\beta_{N-i}}\right]\subset\eS_\e
$$
where $\sigma$ and $D$ are given by the L.S. relation. Note now
that, the $\K_i$, for $i=1,\ldots,n$ normalize $\U^{N}$, and when we
add them to this algebra we perform an iterated  construction of
twisted Laurent polynomial. The resulting algebra will be called
$\T$. We now claim
$$\eS_\e=\T$$
Note that, by construction $\T\subset \eS_\e$, so we only have to
prove that $\eS_\e\subset\T$. Now note that
$$
\E^{k_1}_{\beta_1}\cdots\E^{k_N}_{\beta_N}\K^{s_1}_{1}\cdots\K^{s_n}_{n}\F^{h_N}_{\beta_N}\cdots\F^{h_1}_{\beta_1}\in\T
$$
for every $(k_1,\ldots,k_N)$, $(h_1,\ldots,h_N)\in(\Z^{+})^{N}$ and $(s_1,\ldots,s_n)\in\Z^n$. Then by proposition \ref{PBWbasisS} we have $\eS_\e\subset\T$.
\end{proof}

We finish this section with some remarks on the center of $\U^t_\e$.
Recall that $\U^{t}_{\e}$ is isomorphic to $\U_{\e}$ for every
$t\in\C^*$, hence $Z^t_\e$ is isomorphic to $Z^1_\e=Z_\e$. For
$t=0$, we define $C_{0}$ the subalgebra of $\eS_\e$ generated by
$\E^{l}_{\beta}$, $\F^{l}_{\beta}$ for $\beta\in R^+$ and $\K^{\pm
l}_j$ for $j=1,\ldots,n$ and let $C_\e$ be the center of $\eS_\e$.
Let $Z_0[t]$ the trivial deformation of $Z_0$
\begin{lemma}\label{C0finito}
\begin{enumerate}
    \item $\rho:Z_{0}[t]\rightarrow \U^t_{\e}$ defined in the obvious way is an injective homomorphism of algebra.
    \item $\U^t_\e$ is a free $Z_{0}[t]$ module with base the set of monomials
    $$
    E^{k_1}_{\beta_1}\cdots E^{k_N}_{\beta_N}K^{s_1}_{1}\cdots K^{s_n}_nF^{h_N}_{\beta_N}\cdots F^{h_1}_{\beta_1}
    $$
    for which $0\leq k_i,s_j,h_i<l$, for $i=1,\ldots,N$ and $j=1,\ldots,n$.
\end{enumerate}
\end{lemma}
\begin{proof}
(i) follows by definitions of $Z_0[t]$. (ii) follows from the PBW
theorem.
\end{proof}
\begin{lemma}\label{isocenter}
\begin{enumerate}
\item $Z_0\cong C_{0}$.
\item $\U_\e$ and $\eS_\e$ are isomorphic has $Z_0$ modules.
\end{enumerate}
\end{lemma}
\begin{proof}
Follows from the definitions.
\end{proof}

\subsection{General case}

We can now study the general case.

\begin{definition}
Let $\U^t_\e(\p)$ be the subalgebra of $\U^t_\e$ generated by
$E_{\beta^1_i}$, $F_{\beta_j}$ and $K^{\pm1}_{s}$ for
$i=1,\ldots,h$, $j=1,\dots,N$ and $s=1,\ldots,n$.
\end{definition}
Set $\eS_\e(\p)=\U^{t=0}_{\e}(\p)\subset\eS_\e$.
\begin{proposition}
\begin{enumerate}
    \item For every $t\in\C$, $\U^t_\e(\p)$ is a Hopf subalgebra of $\U^t_\e$.
    \item For any $\la\neq 0$, $\vartheta_\la$ defines by \ref{vartheta} is an algebra isomorphism between $\U^t_\e(\p)$ and $\U^{\la t}_{\e}(\p)$.
\end{enumerate}
\end{proposition}
\begin{proof}
This is an immediate consequence of the same properties in the case
$\p=\g$.
\end{proof}
We can now state the main theorem of this section
\begin{theorem}\label{thmdefpara}
$\eS_\e(\p)$ is a twisted derivation algebra
\end{theorem}
\begin{proof}
We use the same technique as we used in the proof of theorem
\ref{defor}. Let $\D(\p)=\U_\e(\bi^{\el}_+)\otimes \U_\e(\bi_-)$. Define
$$
\Sigma:\eS_\e(\p)\rightarrow\D(\p)
$$
by $\Sigma(E_i)=\E_i$, $\Sigma(F_j)=\F_j$, $\Sigma(K^{\pm1}_j)=\K^{\pm1}_j$ for $i\in\Pi^\el$ and $j=1,\ldots,n$.
\begin{lemma}
$\eS_\e(\p)$ is a subalgebra of $\eS_\e$
\end{lemma}
\begin{proof}
Note that $\D(\p)$ is a subalgebra of $\D$, and, as in lemma
\ref{sigma}, the map $\Sigma$ is well defined and injective. So, we
have the following commutative diagram
$$
\begin{diagram}
\node{\eS_\e(\p)}\arrow{e,t}{\Sigma}\arrow{s,l}{i}\node{\D(\p)}\arrow{s,r}{j}\\
\node{\eS_\e}\arrow{e,t}{\Sigma}\node{\D}
\end{diagram}
$$
Since $\Sigma$ and $j$ are injective maps, we have that $i$ is also
injective
\end{proof}
So we can identify $\eS_\e(\p)$ with the subalgebra of $\eS_\e$
generated by $\E_{\beta^1_i}$, $\F_{\beta_s}$ and $\K^{\pm1}_{j}$
for $i=1,\ldots,h$, $s=1,\ldots,N$ and $j=1,\ldots,n$. As a corollary
of proposition \ref{PBWbasisS} and LS relations, we
have:
\begin{proposition}\label{PBWbasisSp}
\begin{enumerate}
    \item The monomials
$$
\E^{k_1}_{\beta^1_1}\ldots \E^{k_h}_{\beta^1_h}\K^{s_1}_{1}\ldots \K^{s_n}_{n}\F^{t_1}_{\beta_N}\ldots \F^{t_1}_{\beta_1}
$$
for $(k_1,\ldots,k_h)\in(\Z^{+})^{h}$, $(t_1,\ldots,t_N)\in(\Z^{+})^{N}$ and $(s_1,\ldots,s_n)\in\Z^{n}$, form a $\C$ basis of $\eS_\e(\p)$.
    \item For $i<j$ one has
\begin{enumerate}
    \item
    \begin{equation}
    \E_{\beta^1_j}\E_{\beta^1_i}-\e^{(\beta^1_i|\beta^1_j)}\E_{\beta^1_i}\E_{\beta^1_j}=\sum_{k\in\Z^{N}_{+}}c_k\E^k
    \end{equation}
    where $c_k\in \C$ and $c_k\neq 0$ only when $k=\left(k_1,\ldots,k_h\right)$ is such that $k_s=0$ for $s\leq i$ and $s\geq j$, and $\E^k=\E_{\beta^1_1}^{k_1}\ldots \E_{\beta^1_h}^{k_h}$.
    \item
    \begin{equation}
    \F_{\beta_j}\F_{\beta_i}-\e^{-(\beta_i|\beta_j)}\F_{\beta_i}\F_{\beta_j}=\sum_{k\in\Z^{N}_{+}}c_k\F^k
    \end{equation}
    where $c_k\in \C$ and $c_k\neq 0$ only when $k=\left(k_1,\ldots,k_N\right)$ is such that $k_s=0$ for $s\leq i$ and $s\geq j$, and $\F^k=\F_{\beta_N}^{k_N}\ldots \F_{\beta_1}^{k_1}$.
\end{enumerate}
\end{enumerate}
\end{proposition}
So we have:
\begin{theorem}\label{flat}
The monomials
$$
\E^{k_1}_{\beta^1_1}\ldots \E^{k_h}_{\beta^1_h}\K^{s_1}_{1}\ldots \K^{s_n}_{n}\F^{t_N}_{\beta_N}\ldots \F^{t_1}_{\beta_1}
$$
for $(k_1,\ldots,k_h)\in(\Z^{+})^{h}$,
$(t_1,\ldots,t_N)\in(\Z^{+})^{N}$ and $(s_1,\ldots,s_n)\in\Z^{n}$,
are a $\C[t]$ basis of $\U^t_\e$. In particular $t$ is not a zero
divisor in $\U^t_\e$ hence $\U^t_\e$ is a flat over $\C[t]$
\end{theorem}

As we see in general case, we can conclude that $\eS_\e(\p)$ is a quasi polynomial algebra.
\end{proof}

\section{The degree}

\subsection{The degree of $\eS_\e(\p)$}

Using the method exposed in \cite{dp} we can now start the calculation of the degree of $\U_\e(\p)$. Let $\theta=\sum^n_{i=1}a_i\alpha_j$ be the longest root we say that $l$ is a good integer if $l$ is coprime with $a_i$ for all $i$, $2$ and $l$ is coprime also with $6$ if $\g$ has some factor of type $G_2$. So we can now state the main theorem of this section.

\begin{theorem}\label{degS}
If $l$ is a good integer, then
$$
\deg\eS_\e(\p)=l^{\frac{1}{2}\left(l\left(w_0\right)+l\left(w^\el_0\right)+\rk\left(w_0-w^\el_0\right)\right)}
$$
\end{theorem}
\begin{proof}
Denote by $\overline{\eS}_\e(\p)$ the quasi polynomial algebra
associated to $\eS_\e(\p)$. We know by the general theory that
$$
\deg\eS_\e(\p)=\deg\overline{\eS}_\e(\p).
$$
Let $x_i$ denote the class of $E_{\beta^1_i}$ in $\overline{\eS}_\e(\p)$ for $i=1,\ldots,h$ and $y_j$ the class of $F_{\beta_j}$ for $j=1,\ldots,N$, then from theorem \ref{thmdefpara} we have
\begin{eqnarray}
    x_ix_j=\e^{(\beta^1_i|\beta^1_j)}x_jx_i,\\
    y_iy_j=\e^{-(\beta_i|\beta_j)}y_jy_i.
\end{eqnarray}
if $i<j$. Thus we introduce the skew symmetric matrices
$A=\left(a_{ij}\right)$ with $a_{ij}=(\beta_i|\beta_j)$ for $i<j$
and $A^\el=\left(a'_{ij}\right)$ with
$a'_{ij}=(\beta^1_i|\beta^1_j)$ for $i<j$.

Let $k_i$ be the class of $K_i$, using the relation in theorem
\ref{thmdefpara} we obtain a $n\times N$ matrix
$B=\left((w_i|\beta_j)\right)$ and a $h\times N$ matrix
$B^{\el}=((w_i|\beta^1_j))$.

Let $t=2$ unless the Cartan matrix is of type $G_2$, in which case
$t=6$. Since we will eventually reduce modulo $l$ an odd integer
coprime with $t$, we start inverting $t$. Thus consider the free
$\Z\left[\frac{1}{t}\right]$ module $V^+$ with basis
$u_1,\ldots,u_h$, $V^-$ with basis $u'_1,\ldots,u'_N$ and $V^0$ with
basis $w_1,\ldots,w_n$. On $V=V^+\oplus V^0\oplus V^-$ consider the
bilinear form given by
$$
T=\left(\begin{array}{ccc}
A^\el & -^{t}B^\el & 0\\
B^\el & 0& -B\\
0 & ^{t}B & -A
\end{array}\right),
$$
then the rank of $T$ is the degree of $\overline{\eS}_\e(\p)$.

Consider the operators $M^\el=\left(\begin{array}{ccc}
A^\el & -^{t}B^\el & 0
\end{array}\right)$, $M=\left(\begin{array}{ccc}
0 & ^{t}B & -A
\end{array}\right)$, and $N=\left(\begin{array}{ccc}
B^\el & 0& -B
\end{array}\right)$, so that $T=M^\el\oplus N\oplus M$.

Note that
$$
B(u'_i)=\beta_i
$$
and
$$
B^\el(u_i)=\beta^1_i
$$

Now we need some technical lemma:

\begin{lemma}\label{wdeco}
Let $w\in\W$ and fix a reduced expression $w=s_{i_1}\cdots s_{i_k}$.
Given $\omega=\sum^{n}_{i=1}\delta_i\w_i$, with $\delta_i=0$ or $1$.
Set
$$
I_\omega(w):=\left\{t\in\left\{1,\ldots,k\right\}:s_{i_t}(\omega)\neq\omega\right\}.
$$
Then
$$
\omega-w(\omega)=\sum_{t\in I_\omega}\beta_t
$$
\end{lemma}
\begin{proof}
We proceed by induction on the length of $w$. The hypothesis made
implies $s_i(\omega)=\omega$ or $s_i(\omega)=\omega-\alpha_i$. Write
$w=w's_{i_k}$. If $k\notin I_\omega$, then $w(\omega)=w'(\omega)$
and we are done by induction. Otherwise
$$
w(\omega)=w'(\omega-\alpha_{i_k})=w'(\omega)-\beta_k
$$
and again we are done by induction.
\end{proof}

\begin{lemma}\label{directsum}
Let $\theta=\sum^{n}_{i=1}a_i\alpha_i$ the highest root of the root
system $R$. Let $\Z''=\Z[a^{-1}_{1},\ldots,a^{-1}_{n}]$, and let
$\Lambda''=\Lambda\otimes_{\Z}\Z''$ and $Q''=Q\otimes_{\Z}\Z''$.
Then the $\Z''$ submodule $(1-w_{0})\Lambda''$ of $Q''$ is a direct
summand.
\end{lemma}
\begin{proof}
See \cite{dkp2} or \cite{dp} $§10$.
\end{proof}

\begin{lemma}\label{kerM}
\begin{enumerate}
    \item\label{Tsurj} The operator $M$ is onto $V^0\oplus V^-$
    \item The vector $v_\omega:=\left(\sum_{t\in I_\omega}u_t\right)-\omega-w_0(\omega)$, as $\omega$ run thought the fundamental weights, form a basis of the kernel of $M$.
    \item $N(v_\omega)=\omega-w_0(\omega)=\sum_{t\in I_\omega}\beta_t$.
\end{enumerate}
\end{lemma}
\begin{proof}
See \cite{dkp2} or \cite{dp} $§10$.
\end{proof}

Notices that the same result is worth for $M^\el$. Set $T_1=M^\el\oplus M$, then using the notation of lemma \ref{wdeco}, we have
\begin{lemma}\label{kerT1}
The vector $v_\w=\sum_{t\in I_\w(w^\el_0)}u_t-\w-w_0(\w)+\sum_{t\in
I_\w(w_0)}u'_t$, as $\w$ runs through the fundamental weights, form
a basis of the kernel of $T_1$.
\end{lemma}
\begin{proof}
First, we observe that $T_1$ is onto, since $M$ and $M^\el$
are projections over $V^-$ and $V^+$ respectively, by lemma
\ref{kerM}. Since the $n$ vectors $v_\w$ are part of a basis and,
the kernel of $T_1$ is a direct summand of rank $n$, by
surjectivity. It is enough to show that $v_\w$ is in the kernel of
$T_1$. We have
\begin{eqnarray*}
T_1(v_\w)&=&A^\el\left(\sum_{t\in I_\w(w^\el_0)}u_t\right)-^tB^\el\left(-\w-w_0(\w)\right)\\
&&+^tB\left(-\w-w_0(\w)\right)-A\left(\sum_{t\in I_\w(w_0)}u'_t\right)\\
&=&M^\el\left(v_\w\right)-^tB^\el\left(w^\el_0(\w)-w_0(\w)\right)-M(v_\w).
\end{eqnarray*}
So from lemma \ref{wdeco} and lemma \ref{kerM}, we have:
$$
T_1(v_w)=-^tB^\el\left(w^\el_0(\w)-w_0(\w)\right)
$$
Let $w_0=w^\el_0\overline{w}$, since $w$ runs through the
fundamental weights, we have two cases:
\begin{enumerate}
    \item $\overline{w}(\w)=w$, therefore $w^\el_0(\w)-w_0(\w)=0$ and $T_1(v_\w)=0$.
    \item $\overline{w}(\w)\neq \w$, therefore $w^\el_0(\w)=\w$ and $w^\el_0(\w)-w_0(\w)=\w-w_0(\w)\in\ker {^tB^\el}$, by definition of $^tB^\el$, so $T_1(v_\w)=0$.
\end{enumerate}

\end{proof}
Since $T$ is the direct sum of $T_1$ and $N$, its kernel is the
intersection of the $2$ kernels of these operators. We have computed
the kernel of $T_1$ in lemma \ref{kerT1}. Thus the kernel of $T$ equals
the kernel of $N$ restricted to the submodule spanned by the $v_\w$.
\begin{lemma}
$$
N(v_\w)=\sum_{t\in I_\w(w^\el_0)}\beta^1_t-\sum_{t\in I_\w(w_0)}\beta_t=w_0(\w)-w^{\el}_{0}(\w).
$$
\end{lemma}
\begin{proof}
Note that $B(u_t)=\beta_t$, then
\begin{eqnarray}
    N(v_w)=\sum_{t\in I_\w(w^\el_0)}\beta^1_t-\sum_{t\in I_\w(w_0)}\beta_t.
\end{eqnarray}
Finally, the claim follows using lemma \ref{wdeco}.
\end{proof}
Thus, we can identify $N$ we the map
$w_0-w^{\el}_{0}:P\rightarrow Q$. At this point we need the
following fact
\begin{lemma}
Let $\theta=\sum^{n}_{i=1}a_i\alpha_i$ the highest root of the root
system $R$. Let $\Z'=\Z[a^{-1}_{1},\ldots,a^{-1}_{n}]$, and let
$P'=P\otimes_{\Z}\Z'$ and $Q'=Q\otimes_{\Z}\Z'$. Then
the $\Z'$ submodule $(w_0-w^{\el}_{0})P'$ of $Q'$ is a direct
summand.
\end{lemma}
\begin{proof}
The claim follows as a consequence of lemma \ref{directsum}.
\end{proof}
So if $l$ is a good integer, i.e. $l$ is coprime with $t$ and $a_i$
for all $i$, we have
$$
\rk T=l(w_0)+l(w^{\el}_{0})+n-\left(n-\rk\left(w_0-w^{\el}_{0}\right)\right),
$$
and so the theorem follows.
\end{proof}

\subsection{The degree of $\U_\e(\p)$}\label{mainresult}

We begin by observing that every irreducible $\U_\e(\p)$ module $V$
is finite dimensional. Indeed, let $\call{Z}(V)$ be the subalgebra
of the algebra of intertwining operators of $V$ generated by the
action of the elements in $Z_\e(\p)$. Since $\U_\e(\p)$ is finitely
generated as a $Z_\e(\p)$ module, $V$ is finitely generated as
$\call{Z}(V)$ module. If $0\neq f\in\call{Z}(V)$, then $f\cdot V=V$,
otherwise $f\cdot V$ is a proper submodule $V$. Hence, by Nakayama's
lemma, there exist an endomorphism $g\in\call{Z}(V)$ such that
$1-gf=0$, i.e. $f$ is invertible. Thus $\call{Z}(V)$ is a field. It
follows easily that $\call{Z}(V)$ consists of scalar operators. Thus
$V$ is a finite dimensional vector space.

Since $Z_\e(\p)$ acts by scalar operators on $V$, there exists an
homomorphism $\chi_V:Z_\e\mapsto\C$, the \emph{central character of
V}, such that
$$
z\cdot v=\chi_V(z)v
$$
for all $z\in Z_\e$ and $v\in V$. Note that isomorphic
representations have the same central character, so assigning to a
$\U_\e(\p)$ module its central character gives a well defined map
$$
\Xi:\Rap\left(\U_\e(\p)\right)\rightarrow\Spec\left(Z_\e(\p)\right),
$$
where $\Rap\left(\U_\e(\p)\right)$ is the set of isomorphism classes
of irreducible $\U_\e(\p)$ modules, and $\Spec\left(Z_\e(\p)\right)$
is the set of algebraic homomorphisms $Z_\e(\p)\mapsto\C$.

To see that $\Xi$ is surjective, let $I^{\chi}$, for $\chi\in\Spec\left(Z_\e(\p)\right)$, be the ideal in $\U_\e(\p)$ generated by
$$
\ker\chi=\left\{z-\chi(z)\cdot 1:z\in Z_\e(\p)\right\}.
$$
To construct $V\in\Xi^{-1}(\chi)$ is the same as to construct an
irreducible representation of the algebra
$\U^{\chi}_{\e}(\p)=\U_\e(\p)/I^\chi$. Note that
$\U^{\chi}_{\e}(\p)$ is finite dimensional and non zero. Thus, we
may take $V$, for example, to be any irreducible subrepresentation
of the regular representation of $\U^{\chi}_{\e}(\p)$.

Let $\chi\in\Spec(Z_0(\p))$, we define,
$$
\U^{\chi}_{\e}(\p)=\U_{\e}(\p)/J^\chi
$$
where $J^\chi$ is the two sided ideal generated by
$$
\ker\chi=\left\{z-\chi(z)\cdot 1:z\in Z_0(\p)\right\}
$$

As we have seen at the end of section \ref{gcaso},
$Z_0(\p)[t]\subset\U^t_\e(\p)$, so for all $t\in\C$ and $\chi\in
\Spec(Z_{0}(\p))$, we can define
$\U^{t,\chi}_{\e}(\p)=\U^{t}_{\e}(\p)/J^\chi$ where $J^\chi$ is the
two side ideal generated by
$$
\ker\chi=\left\{z-\chi(z)\cdot 1:z\in Z_0(\p)\right\}
$$

The PBW theorem for $\U^{t}_{\e}(\p)$ implies that
\begin{proposition}
The monomials
$$
    E^{s_{1}}_{\beta^{1}_{1}}\cdots E^{s_{h}}_{\beta^{1}_{h}}K^{r_1}_{1}\cdots K^{r_n}_{n}F^{t_{k+h}}_{\beta^{1}_{h}}\cdots F^{t_{k+1}}_{\beta^{1}_{1}}F^{t_k}_{\beta^{2}_{k}}\cdots F^{t_1}_{\beta^{2}_{1}}
    $$
    for which $0\leq s_j,t_i,r_v<l$, form a $\C$ basis for
$\U^{t,\chi}_{\e}(\p)$
\end{proposition}

\begin{lemma}
For every $0\neq\la\in\C$, $\U^{t,\chi}_{\e}(\p)$ is isomorphic to
$\U^{\la t,\chi}_{\e}(\p)$.
\end{lemma}
\begin{proof}
Consider the isomorphism $\vartheta_\la$ from $\U^{t}_{\e}(\p)$ and
$\U^{\la t}_{\e}(\p)$, defined by (\ref{vartheta}). Its
follows from the above definition that
$\vartheta_\la(J^{\chi})=J^{\chi}$. Then $\vartheta_\la$ induce an
isomorphism between $\U^{t,\chi}_{\e}$ and $\U^{\la t,\chi}_{\e}$.
\end{proof}

\begin{proposition}
The $\U_{\e}(\p)$ algebras $\U^{t,\chi}_{\e}(\p)$ form a continuous
family parametrized by $\Ze=\C\times \Spec(Z_0(\p))$.
\end{proposition}
\begin{proof}
Let $\V$ denote the set of triple $(t,\chi,u)$ with $(t,\chi)\in\Ze$
and $u\in\U^{t,\chi}_{\e}(\p)$. Then from the PBW theorem we have
that the set of monomial
$$
    E^{s_1}_{\beta_1}\cdots E^{s_h}_{\beta^1_h}K^{r_1}_{1}\ldots K^{r_n}_{n} F^{t_N}_{\beta_N}\cdots F^{t_1}_{\beta_1}
$$
for which $0\leq s_i, t_i, r_v<l$, for $i\in\Pi^\el$, $j=1,\ldots,N$
and $v=\ldots,n$, form a basis for each algebra $\U^{t,\chi}_{\e}$.

Order the previous monomials and assign to $u\in\U^{t,\chi}_{\e}$
the coordinate vector of $u$ with respect to the ordered basis. This
construction identifies $\V$ with $\Ze\times\C^d$, where
$d=l^{h+n+N}$, thereby giving $\A$ a structure of an affine variety.

Consider the vector bundle $\pi:\V\rightarrow \Ze$,
$(t,\chi,u)\rightarrow(t,\chi)$. Note that the structure constant of
the algebra $\U^{t,\chi}_{\e}(\p)$, as well as the matrix entries of
the linear transformations which define the action of $\U_\e(\p)$
relative to the basis, are polynomial in $\chi$ and $t$. This means
that the maps
\begin{eqnarray*}
    &\mu:\V\times_\Ze\V\rightarrow\V,\mbox{ }&\left((t,\chi,u),(t,\chi,v)\right)\mapsto(t,\chi,uv)\\
    &\rho:\U_{\e}\times\V\rightarrow\V,\mbox{   }&\left(x,(t,\chi,u)\right)\mapsto(t,\chi,x\cdot u)
\end{eqnarray*}
where $(t,\chi)\in\Ze$, $u$, $v\in\U^{t,\chi}_{\e}(\p)$ and
$x\in\U_{\e}(\p)$, define on $\V$ a structure of vector
bundle of algebra and a structure of vector bundle of
$\U_\e(\p)$-modules. The fiber of $\pi$ above $(t,\chi)$ is the
$\U_\e(\p)$-algebra $\U^{t,\chi}_{\e}(\p)$.
\end{proof}

If we fix $\chi\in\Spec(Z_0)$, we have from theorem \ref{flat} that
the family of algebra $\U^{t,\chi}_{\e}(\p)$ is a flat deformation
of algebra over $\Spec\C[t]$.

Summarizing, if $\e$ is a primitive $l^{th}$ root of $1$ with $l$
odd and $l>d_i$ for all $i$, we have prove the following facts on
$\U_\e$:
\begin{itemize}
    \item $\U^t_\e$ and $U_\e(\p)$ are domains because $\U_\e(\g)$ it is,
    \item $\U^t_\e$ and $\U_\e(\p)$ are finite modules over $Z_0[t]$ and $Z_0$ respectively (cf lemma \ref{C0finito} and proposition \ref{center2}).
\end{itemize}
Since the LS relations holds for $\U_\e(\p)$ and $\U^t_\e$ (cf
proposition \ref{LSp}), we can apply the theory developed in \cite{dp}, and we obtain that $\Gr\U_\e(\p)$ and $\Gr\U^t_\e$ are
twisted polynomial algebra, with some elements inverted. Hence all
conditions of the characterization of maximal order (see theorem 6.4 in \cite{dp} ) are verified, we have
\begin{theorem}\label{maxordUp}
$\U^t_\e$ and $\U_\e(\p)$ are maximal orders.
\end{theorem}
Therefore, $\U_\e(\p)\in\Ci_m$, i.e. is an algebra with trace of
degree $m$.
\begin{theorem}\label{thmdegU}
The set
\begin{eqnarray*}
    \Omega=\left\{a\in \Spec(Z_\e(\p)),\mbox{ such that the corresponding semisimple}\right.\\
    \left.\mbox{representation of $\U_\e(\p)$ is irreducible}\right\}
\end{eqnarray*}
is a Zariski open set. This is exactly the part of $\Spec(Z_\e(\p))$
over which $\U^{\chi}_{\e}(\p)$ is an Azumaya algebra of degree $m$.
\end{theorem}
\begin{proof}
Apply theorem 4.5 in \cite{dp}, with $R=\U_\e(\p)$ and $T=Z_\e(\p)$.
\end{proof}

Recall that $Z_\e(\p)$ is a finitely generated module over
$Z_0(\p)$. Let $\tau:\Spec(Z_\e(\p))\rightarrow\Spec(Z_0(\p))$ be the
finite surjective morphism induced by the inclusion of $Z_0(\p)$ in
$Z_\e(\p)$. The properness of $\tau$ implies the following
\begin{corollary}
The set
$$
\Omega_0=\left\{a\in\Spec(Z_0(\p)):\tau^{-1}(a)\subset\Omega\right\}
$$
is a Zariski dense open subset of $\Spec(Z_0)$.
\end{corollary}

We know by the theory developed in \cite{dp}, that $\eS_\e(\p)\in\Ci_{m_0}$, with $m_0=l^{l(w_0)+l(w^{\el}_{0})+\rk\left(w_0-w^{\el}_{0}\right)}$. As
we see in proposition \ref{C0finito}, $S_\e(\p)$ is a finite module
over $C_0$, then $C_{\e}$, the center of $\eS_\e(\p)$ is finite over
$C_0$. The inclusion $C_0\hookrightarrow C_\e$ induces a projection
$\upsilon:\Spec(C_\e)\rightarrow\Spec(C_0)$. As before, we have:
\begin{lemma}\label{thmdeg}
\begin{enumerate}
\item
\begin{eqnarray*}
    \Omega'=\left\{a\in \Spec(C_\e),\mbox{such that the corresponding semisimple}\right.\\
    \left.\mbox{representation of $\eS_\e(\p)$ is irreducible}\right\}
\end{eqnarray*}
is a Zariski open set. This is exactly the part of $\Spec(C_\e)$ over which $\eS^{\chi}_{\e}(\p)$ is an Azumaya algebra of degree $m_0$.
\item The set
$$
\Omega'_0=\left\{a\in\Spec(Z_0(\p)):\upsilon^{-1}(a)\subset\Omega'\right\}
$$
is a Zariski dense open subset of $\Spec(Z_0)$.
\end{enumerate}
\end{lemma}
\begin{proof}
Apply theorem 4.5 in \cite{dp} at $\eS_{\e}(\p)$.
\end{proof}

Since $\Spec(Z_0)$ is irreducible, we have that
$\Omega_0\cap\Omega'_0$ is non empty.

We can state the main theorem of this section
\begin{theorem}
If $l$ is a good integer, then
$$
\deg\U_{\e}(\p)=l^{\frac{1}{2}\left(l(w_0)+l(w^{\el}_{0})+\rk\left(w_0-w^{\el}_{0}\right)\right)}
$$
\end{theorem}
\begin{proof}
For $\chi\in\Omega_0\cap\Omega'_0$, we have, using theorem
\ref{thmdegU} and lemma \ref{thmdeg},
$$
\deg\U_{\e}(\p)=m=\deg\U^{\chi}_\e(\p),
$$
and
$$
\deg\eS^\chi_\e(\p)=\deg\eS_\e(\p).
$$
But for all $t\neq0$, we have that $\U^{t,\chi}_{\e}$ is isomorphic
to $\U^{\chi}_{\e}(\p)$ as algebra. By construction $\U^{\chi}_\e(\p)$ as a module over itself is irreducible, hence it is a simple algebra. By rigidity of semisimple algebra (\cite{pierce} or \cite{procesi3}) we have that $\eS^\chi_\e(\p)=\U^{0,\chi}_{\e}$ is isomorphic to $\U^{\chi}_{\e}(\p)$. Then
$$
\deg\U_{\e}(\p)=m=\deg\U^{\chi}_\e(\p)=\deg\eS^\chi_\e(\p)=\deg\eS_\e(\p).
$$
And by theorem \ref{degS} the claim follows.
\end{proof}

As $\U_\e(\p)$ is a maximal order, $Z_\e(\p)$ is integrally closed, so we can make the following construction: denote by
$Q_\e:=Q(Z_\e(\p))$ the field of fractions of $Z_\e(\p)$, we have
that $Q(\U_\e(\p))=\U_\e(\p)\otimes_{Z_\e(\p)}Q_\e$ is a division
algebra, finite dimensional over its center $Q_\e$. Denote by $\F$
the maximal commutative subfield of $Q(\U_\e(\p)$, we have, using
standard tools of associative algebra (cf \cite{pierce}), that
\begin{enumerate}
    \item $\F$ is a finite extension of $Q_\e$ of degree $m$,
    \item $Q(\U_\e(\p))$ has dimension $m^2$ over $Q_\e$,
    \item $Q(\U_\e(\p))\otimes_{Q_\e}\F\cong M_m(\F)$.
\end{enumerate}
Hence, we have that
\begin{eqnarray*}
    \dim_{Q(Z_0(\p))}(Q_\e)&=&\deg(\tau)\\
    \dim_{Q_\e}(Q(\U_\e(\p)))&=&m^2\\
    \dim_{Q(Z_0(\p))}(Q(\U_\e(\p)))&=&l^{h+N+n}
\end{eqnarray*}
where, the first equality is a definition, the second has been pointed out above and the third follows from the P.B.W theorem. Then, we have
$$
l^{h+N+n}=m^2\deg(\tau)
$$
with $m=l^{\frac{1}{2}\left(l(w_0)+l(w^{\el}_{0})+\rk\left(w_0-w^{\el}_{0}\right)\right)}$, so
\begin{corollary}
$$
\deg(\tau)=l^{n-\rk\left(w_0-w^{\el}_{0}\right)}.
$$
\end{corollary}

\section{The center}

We want to make some observations on the center of $\U^{\chi,t}_\e$ and $\U^t_\e$ that perhaps they can be useful in the explicit determination of the center of $\U_\e(\g)$. 

\subsection{The center of $\U^\chi_\e(\p)$}

We want to explain a method which in principle allows us to
determine the center of $\U^\chi_\e(\p)$ for all $\chi\in\Spec(Z_0)$.

Let $\chi_0\in\Spec(Z_0(\p))$ define by $\chi_0(E_i)=0$, $\chi_0(F_i)=0$ and $\chi_0(K^{\pm1}_{i})=1$, we set $\U^{0}_{\e}(\p)=\U^{\chi_0}_{\e}(\p)$.
\begin{proposition}
$\U^{0}_{\e}(\p)$ is a Hopf algebra with the comultiplication, counit and antipode induced by $\U_{\e}(\p)$.
\end{proposition}
\begin{proof}
This is immediately since $J^{\chi_0}$ is an Hopf ideal.
\end{proof}

\begin{proposition}
let $\chi\in\Spec(Z_0(\p))$
\begin{enumerate}
    \item $\U^{\chi}_{\e}(\p)$ is an $\U_{\e}(\p)$ module, with the action define by
    $$
    a\cdot u=\sum a_{(1)}uS(a_{(2)}),
    $$
    where $\Delta(a)=\sum a_{(1)}\otimes a_{(2)}$.
    \item $\U^{\chi}_{\e}(\p)$ is an $\U^{0}_{\e}(\p)$ module, with the action induced by $\U_{\e}(\p)$.
\end{enumerate}
\end{proposition}
\begin{proof}
Easy verification of the proprieties.
\end{proof}

\begin{proposition}
Let $x\in\U^{\chi}_{\e}(\p)$. Then $x$ is in the center of $\U^{\chi}_{\e}(\p)$ if and only if $x$ is invariant under the action of $\U^{0}_{\e}(\p)$, that is
\begin{eqnarray}\label{invariantp}
    E_i\cdot x&=&0,\\
    F_i\cdot x&=&0,\\
    K_i\cdot x&=&x.
\end{eqnarray}
\end{proposition}
\begin{proof}
Let $x\in Z\left(\U^{\chi}_{\e}(\p)\right)$ then
$$
E_i\cdot x=E_ix-K_ixK^{-1}_{i}E_i=0,
$$
in the same way we obtain the other relations.

Suppose now that $x$ verify the relations \ref{invariantp}. Then
$$
K_i\cdot x=K_ixK^{-1}_{i}=x,
$$
imply that $K_ix=xK_i$. From $E_i\cdot x=0$ we obtain
\begin{eqnarray*}
0=E_i\cdot x&=&E_ix-K_ixK^{-1}_{i}E_i\\
&=&E_ix-xE_i.
\end{eqnarray*}
its follows that $E_ix=xE_i$. In the same way we have $F_ix=xF_i$. Then $x$ lies in the center.
\end{proof}
So we can determine the center at $t$ generic by lifting the center of the algebra at $t=0$.

\subsection{The center of $\U^t_\e$}

We want to study  the restriction of the deformation at the center of $\U_\e(\p)$. Since $t$ is without torsion, it is easy to see that,
\begin{proposition} 
\begin{enumerate}
\item $Z_{\e,0}:=Z_\e(\p)[t]/tZ_\e(\p)[t]=Z_\e(\p)[t]/t\U^t_\e\cap Z_\e(\p)[t]$.
\item $Z_{\e,0}\cong Z_\e(\p)$.
\end{enumerate}
\end{proposition}
We want to prove that 
\begin{theorem}
$$
Z_{\e,0}=C_\e(\p).
$$
\end{theorem}

Note that $Z_{\e,0}$ and $C_\e$ are integrally closed domain and finitely generated $Z_0(\p)$ algebras such that:
\begin{lemma}
$$
Q(Z_{\e,0})=Q(C_{\e})
$$
where $Q(Z_{\e,0})$ and $Q(C_{\e})$ are the fields of fractions of $Z_{\e,0}$ and $C_\e$ respectively.
\end{lemma}
\begin{proof}
Note that $Q(Z_{\e,0})\subset Q(C_{\e})$ and $Z_{\e,0}\cong Z_\e(\p)$. Since $\U_\e(\p)$ and $\eS_\e(\p)$ have the same degree, the observation at the and of section \ref{mainresult} implies that
$$
\dim_{Q(Z_0(\p)}Q(Z_{\e,0})=\dim_{Q(Z_0(\p)}Q(C_\e)
$$
and the result follows.
\end{proof}
So using classic results (\cite{serre},\cite{matsumura}) we can conclude that $Z_{\e,0}=C_\e$.

\bibliographystyle{alpha}
\bibliography{qgroups,lie,associative_algebra}

\end{document}